\documentclass[preprint,12pt]{elsarticle}

\usepackage{amssymb}
\usepackage{amsmath}
\usepackage{amsthm}
\newdefinition{dfn}{Definition}
\newdefinition{rmk}{Remark}
\newtheorem{thm}{Theorem}
\newproof{pf}{Proof}
\journal{Journal of Multivariate Analysis}

\begin{document}

\begin{frontmatter}

%% Title, authors and addresses

%% use the tnoteref command within \title for footnotes;
%% use the tnotetext command for the associated footnote;
%% use the fnref command within \author or \address for footnotes;
%% use the fntext command for the associated footnote;
%% use the corref command within \author for corresponding author footnotes;
%% use the cortext command for the associated footnote;
%% use the ead command for the email address,
%% and the form \ead[url] for the home page:
%%
%% \title{Title\tnoteref{label1}}
%% \tnotetext[label1]{}
%% \author{Name\corref{cor1}\fnref{label2}}
%% \ead{email address}
%% \ead[url]{home page}
%% \fntext[label2]{}
%% \cortext[cor1]{}
%% \address{Address\fnref{label3}}
%% \fntext[label3]{}

\title{Multivariate Copula Expressed by Lower Dimensional Copulas}

%% use optional labels to link authors explicitly to addresses:
%% \author[label1,label2]{<author name>}
%% \address[label1]{<address>}
%% \address[label2]{<address>}

\author[avf]{Edith Kov\'{a}cs}

\address[avf]{Department of Mathematics, \'{A}VF College of Management of Budapest,
              Vill\'{a}nyi \'{u}t 11-13, H-1114 Budapest, Hungary\\Corresponding author, \\ E-mail address: kovacs.edith@avf.hu}
              
\author[bme]{Tam\'{a}s Sz\'{a}ntai}

\address[bme]{Institute of Mathematics, Budapest University of Technology and Economics,
              M\H{u}egyetem rkp. 3, H-1111 Budapest, Hungary\\ E-mail address: szantai@math.bme.hu}

\begin{abstract}
%% Text of abstract

Modeling of high order multivariate probability distribution is a
difficult problem which occurs in many fields. Copula approach is a good
choice for this purpose, but the curse of dimensionality still remains a
problem. In this paper we give a theorem which expresses a multivariate
copula by using only some lower dimensional ones based on the conditional
independences between the variables. In general the construction of a
multivariate copula using this theorem is quite difficult, due the
consistency properties which have to be fulfilled. For this purpose we
introduce the sample derivated copula, and prove that the dependence between the
random variables involved depends just on this copula and on the partition.
By using the sample derivated copula the theorem can be successfully applied, in
order to to construct a multivariate discrete copula by using some of its
marginals.

\end{abstract}

\begin{keyword}
%% keywords here, in the form: keyword \sep keyword

%% MSC codes here, in the form: \MSC code \sep code
%% or \MSC[2008] code \sep code (2000 is the default)

Multivariate copula, Junction tree, Conditional independence, Sample derivated copula.

\end{keyword}

\end{frontmatter}

%%
%% Start line numbering here if you want
%%
% \linenumbers

%% main text
\section{Introduction}
\label{sec:1}

First we motivate why should we model the multivariate 
distribution by copulas from an information theoretical point of view.  
The information content of a multivariate probability distribution depends only
on its copula density. In \cite{MaSu08} and \cite{CaVi09} one can see this result for the
two-dimensional case and the same is true for more dimensions, too.

In this paper we prove a theorem which links the multivariate probability
distribution assigned to a junction tree to the multivariate copula. It is
known that the probability distribution assigned to a junction tree uses the
conditional independence structure underlying the random variables so the
copula introduced here will have this property, too.

In this introductory part we describe the main concepts and introduce the
notations which we will use in the paper. In the second section we prove a theorem
which links a multivariate copula to the junction tree probability distribution.
In the third section we will
introduce the concept of Sample Derivated Copula (SDC) which makes possible the 
exploitation of the conditional independences between the random variables. 
We prove that the information content of the probability distribution given 
by a partition set depends only on the SDC. 
In the fourth section we apply the junction tree approach to the SDC.

We finish the paper with conclusions and possible applications.

Let $V=\{1,\ldots,n\}$ be a set of vertices. A hypergraph is a
set $V$ of vertices together with a set $\Gamma$ of subsets of $V$. A
hypergraph is acyclic if no elements in $\Gamma$ are subsets of other
elements, and if the elements of $\Gamma$ can be ordered
$(K_{1},\ldots,K_{m})$ to have the \textit{running intersection property}:
for all $j \geq 2,$ exists $i<j:K_{i}\supseteq K_{j}\cap\left(  K_{1}\cup
\ldots\cup K_{j-1}\right)$ \cite{LaSp84}.

It is convenient to introduce the so called separator sets
$S_{j}=K_{j}\cap(K_{1}\cup\ldots\cup K_{j-1})$, where $S_{1}=\phi$.

We note here that if $R_{j}=K_{j}\backslash S_{j}$ then $S_{j}$
separates (in graph terms) the vertices in $R_{j}$ from the vertices in
$(K_{1}\cup\ldots\cup K_{j-1})\backslash S_{j}$.

We mention here that a hypergraph $(V,\Gamma)$ is acyclic if
and only if $\Gamma$ can be considered to be the set of cliques of a
chordal (triangulated) graph \cite{LaSp88},\cite{Ta84}.

In the following we consider acyclic hypergraphs with the property that the
union of all sets in $\Gamma$ is $V$. We denote the separator set by
$\mathcal{S}$ and refer to the acyclic hypergraph as
$(V,\Gamma,\mathcal{S})$.

Let $V=\{1, 2, \ldots, n\}$ be the set of indices of the continuous random
variables $X=\{X_{1},\ldots,X_{n}\}$. We suppose that the probability density functions
of $X_{1},\ldots,X_{n}$ exist and denote them by $f_{X_{1}},\ldots,f_{X_{n}}$.

We need the following notations:

\begin{itemize}
\item $F_{X_{i}}(x_{i})=P(X_{i}<x_{i};X_{j}=\infty\mbox{ for all }j\neq i)$
stands for the univariate marginal cumulative
distribution function corresponding to the variable $X_{i}$,

\item The joint probability density function and the joint cumulative distribution
function of $(X_{1},\ldots,X_{n})^{T}$ is denoted by
$f_{\mathbf{X}}(\mathbf{x})$ and $F_{\mathbf{X}}(\mathbf{x})$, respectively,

\item $D=\{i_{1},\ldots,i_{d}\} \subset V$,
$\mathbf{X}_{D}=(X_{i_{1}},\ldots,X_{i_{d}})^{T}$,
$\mathbf{x}_{D}=(x_{i_{1}},\ldots,x_{i_{d}})^{T}$,

\item  The $d$-th order marginal probability density function and
the $d$-the order marginal cumulative distribution function
of $\mathbf{X}_{D}$ is denoted by $f_{\mathbf{X}_{D}}(\mathbf{x}_{D})$ and
$F_{\mathbf{X}_{D}}(\mathbf{x}_{D})$,respectively.
\end{itemize}

Having these notations we give the concept of the junction tree. It is known
that the junction tree encodes the conditional independences between the
variables. Let us remark here that from now on the indices of the random
variables are assigned to the nodes of a graph. In the graph a set of nodes $B$
separates a set of nodes $A$ from another set of nodes $C$, where $A,B,C$ are
disjoint subsets of $V$, if and only if $X_{A}$ and $X_{C}$ are conditionally
independent with respect to $X_{B}$ (see the definition of the Markov random field).

\begin{dfn}\label{de:1}
A junction tree over X is a cluster tree, which is assigned to an acyclic
hypergraph $(V,\Gamma,\mathcal{S})$ as follows:
\begin{enumerate}[1)]
\item Each cluster of the cluster tree consists of a subset $X_K$ of $X$ ,
where $K\in\Gamma$. To each cluster is assigned the joint marginal
density function $f_{\mathbf{X}_K}\left( \mathbf{x}_K\right) $;
\item Each edge connecting to clusters is called separator and consists of a
subset $X_S$ of $X$, where $S$ is a separator set. To each separator there
is assigned the marginal probability density function $f_{\mathbf{X}%
_S}\left( \mathbf{x}_S\right) $;
\item The union of all clusters is $X$.
\end{enumerate}
\end{dfn}

\begin{dfn}\label{de:2}
A junction tree probability distribution is a probability distribution
assigned to the junction tree in the following way:
\[
f_{\mathbf{X}}\left( \mathbf{x}\right) =\dfrac{\prod\limits_{K\in\Gamma}f_{%
\mathbf{X}_K}\left( \mathbf{x}_K\right) }{\prod\limits_{S\in\mathcal{S}%
}\left( f_{_{\mathbf{X}_S}}\left( \mathbf{x}_S\right) \right) ^{v_{s-1}}},
\]
where $v_S$ is the number of those clusters which contain all the variables of $%
\mathbf{X}_S$.
\end{dfn}

It is useful to note here that since in the hypergraph $\left(  V,\Gamma
,\mathcal{S}\right)  $ $S_{j}$ separates (in graph terms) the vertices in
$R_{j}=K_{j}-S_{j\mbox{ }}$ from the vertices in $\left(  K_{1}\cup\ldots\cup
K_{j-1}\right)  -S_{j}$ the random variables with indices in $R_{j}%
=K_{j}-S_{j\mbox{ }}$ and the variables with indices in $\left(  K_{1}%
\cup\ldots\cup K_{j-1}\right)  -S_{j}$ are conditionally independent with
respect to the variables with indices in $S_{j}$.

\begin{rmk}\label{re:1}
Since the junction tree is assigned to an acyclic hypergraph, the
running intersection property stands for the junction tree, too. It can be
reformulated as follows. If two clusters contain a random variable, then all
clusters on the path between these clusters contain this random variable.
\end{rmk}

First we call back the concept of copula and formulate the Sklar's theorem
(see \cite{Co89} and \cite{Ne99}).

\begin{dfn}\label{de:3}
A function $C:\left[0;1\right]^{d}\rightarrow\left[0;1\right]$
is called a $d$-dimensional copula if it satisfies the following conditions:

\begin{enumerate}[1)]

\item $C\left(u_{1},\ldots,u_{d}\right)$ is increasing in each component
$u_{i}$,

\item $C\left(u_{1},\ldots,u_{i-1},0,u_{i+1},\ldots,u_{d}\right)=0$ for all
$u_{k}\in\left[0;1\right]$,\ $k\neq i,\ i=1,\ldots,n$,

\item $C\left(  1,\ldots,1,u_{i},1,\ldots,1\right)  =u_{i}$ for all $u_{i}%
\in\left[  0;1\right]  ,\ i=1,\ldots,d$,

\item C is $d$-increasing, i.e for all $\left(  u_{1,1},\ldots,u_{1,d}\right)  $
and $\left(  u_{2,1},\ldots,u_{2,d}\right)  $ in $\left[  0;1\right]  ^{d}$
with $u_{1,i}<u_{2,i}$ for all i, we have

\[
\sum\limits_{i_{1}=1}^{2}\cdots
\sum\limits_{i_{d}=1}^{2}\left(  -1\right)  ^{\sum\limits_{j=1}^{d}i_{j}%
}C\left(  u_{i_{1},1},\ldots,u_{i_{d},d}\right)  \geq0.
\]

\end{enumerate}

\end{dfn}

Due to Sklar's theorem if $X_{1},\ldots,X_{d}$ are continuous random variables
defined on a common probability space, with the univariate marginal cdf's
$F_{X_{i}}\left(  x_{i}\right)  $ and the joint cdf $F_{X_{1},\ldots,X_{d}%
}\left(  x_{1},\ldots,x_{d}\right)  $ then there exists a unique copula
function $C_{X_{1},\ldots,X_{d}}\left(  u_{1},\ldots,u_{d}\right)  :\left[
0;1\right]  ^{d}\rightarrow\left[  0;1\right]  $ such that by the substitution
$u_{i}=F_{i}\left(  x_{i}\right), \ i=1,\ldots,d$ we get

\[
F_{X_{1},\ldots,X_{d}}\left(x_{1},\ldots,x_{d}\right)
=C_{X_{1},\ldots,X_{d}}\left(F_{1}\left(  x_{1}\right)  ,\ldots,F_{d}\left(  x_{d}\right)  \right)
\]
for all $\left(x_{1},\ldots,x_{d}\right)^{T}\in R^{d}.$

In the following we will use the vectorial notation
$F_{\mathbf{X}_{D}}\left(\mathbf{x}_{D}\right)=C_{X_{D}}\left(
\mathbf{u}_{D}\right)$, where $\mathbf{u}_{D}%
=\left(  F_{X_{i_{1}}}\left(  x_{i_{1}}\right)  ,\ldots,F_{X_{i_{d}}}\left(
x_{i_{d}}\right)  \right)  ^{T}$.

We need the following assertion:

\[
\begin{array}{l}
f_{X_{i_{1}},\ldots X_{i_{d}}}\left(  x_{i_{1}},\ldots,x_{i_{d}}\right)=\vspace{3mm}\\
=\dfrac{\partial^{d} F_{X_{i_{1}},\ldots X_{i_{d}}}\left(  x_{i_{1}},\ldots
,x_{i_{d}}\right)  }{\partial x_{i_{1}}\cdots\partial x_{i_{d}}} \vspace{3mm}\\
=\dfrac{\partial^{d} C_{X_{i_{1}},\ldots,X_{i_{d}}}\left(F_{X_{i_{1}}}\left(  x_{i_{1}}\right)  ,\ldots,F_{X_{i_{d}}}\left(  x_{i_{d}}\right)  \right)  }
{\partial x_{i_{1}}\cdots\partial x_{i_{d}}}\vspace{3mm}\\
=\left. \dfrac{\partial^{d} C_{X_{i_{1}},\ldots X_{i_{d}}}\left(  u_{i_{1}},\ldots
,u_{i_{d}}\right)  }{\partial u_{i_{1}}\cdots\partial u_{i_{d}}}
\right|_{\substack{u_{i_{k}}=F_{X_{i_{k}}}\left(  x_{i_{k}}\right), \\ k=1,\ldots
d}}\cdot\prod\limits_{k=1}^{d}\dfrac{\partial F_{X_{i_{k}}}\left(  x_{i_{k}%
}\right)  }{\partial x_{i_{k}}}\vspace{3mm}\\
=c_{X_{i_{1}},\ldots X_{i_{d}}} \left( F_{X_{i_{1}}}\left(x_{i_{1}}\right),\ldots,F_{X_{i_{d}}}\left(  x_{i_{d}}\right)  \right)
\cdot \prod\limits_{k=1}^{d} f_{X_{i_{k}}} \left(x_{i_{k}}\right)
\end{array}
\]
In vectorial notation this can be written as
\begin{equation}\label{eq:1}
f_{\mathbf{X}_{D}}\left(  \mathbf{x}_{D}\right)  =c_{\mathbf{X}_{D}}\left(
\mathbf{u}_{D}\right)  \cdot\prod\limits_{i_{k}\in D}^{{}%
}f_{X_{i_{k}}}\left(  x_{i_{k}}\right)
\end{equation}
and from (\ref{eq:1}) we get

\begin{equation}\label{eq:2}
c_{\mathbf{X}_{D}}\left(  \mathbf{u}_{D}%
\right)  =\dfrac{f_{\mathbf{X}_{D}}\left(  \mathbf{x}_{D}\right)  }%
{\prod\limits_{i_{k}\in D}^{{}}f_{X_{i_{k}}}\left(  x_{i_{k}}\right)  }
\end{equation}

\section{The multivariate copula associated to a junction tree probability distribution.}
\label{sec:2}

\begin{thm}\label{theo:1} 
The copula density function associated to a junction tree
probability distribution
\[
f_{\mathbf{X}}\left(  \mathbf{x}\right)
=\dfrac{\prod\limits_{K\in\Gamma}f_{\mathbf{X}_{K}}\left(  \mathbf{x}%
_{K}\right)  }{\prod\limits_{S\in\mathcal{S}}\left[  f_{\mathbf{X}_{S}}\left(
\mathbf{x}_{S}\right)  \right]  ^{v_{S}-1}},
\]
is given by

\begin{equation}\label{eq:3}
c_{\mathbf{X}}\left(  \mathbf{u}_{V}\right)  = \dfrac{\prod\limits_{K\in\Gamma
}c_{\mathbf{X}_{K}}\left(  \mathbf{u}_{K}\right)  }{\prod\limits_{S\in
\mathcal{S}}\left[  c_{\mathbf{X}_{S}}\left(  \mathbf{u}_{S}\right)  \right]
^{v_{S}-1}}.
\end{equation}
\end{thm}

\begin{pf}

\begin{equation}\label{eq:4}
f_{\mathbf{X}}\left(  \mathbf{x}\right)  =\dfrac{\prod\limits_{K\in\Gamma
}f_{\mathbf{X}_{K}}\left(  \mathbf{x}_{K}\right)  }{\prod\limits_{S\in
\mathcal{S}}\left[  f_{\mathbf{X}_{S}}\left(  \mathbf{x}_{S}\right)  \right]
^{v_{S}-1}}=\dfrac{\prod\limits_{K\in\Gamma}c_{\mathbf{X}_{K}}\left(
\mathbf{u}_{K}\right)  \cdot\prod\limits_{i_{k}\in K}^{{}%
}f_{X_{i_{k}}}\left(  x_{i_{k}}\right)  }{\prod\limits_{S\in\mathcal{S}%
}\left[  c_{\mathbf{X}_{S}}\left(  \mathbf{u}_{S}\right)
\cdot\prod\limits_{i_{k}\in S}^{{}}f_{X_{i_{k}}}\left(  x_{i_{k}}\right)
\right]  ^{v_{S}-1}}.
\end{equation}

The question that we have to answer is how many times appears in the
nominator respectively in the denominator the probability density function
$f_{X_{i}}\left(  x_{i}\right)$ of each $X_{i}$ random variable.

Since $\bigcup\limits_{K\in\Gamma}\mathbf{X}_{K}=X$ for each random variable
$X_{i}$ $\epsilon X$, $f_{X_{i}}\left(  x_{i}\right)$ appears at least once
in the nominator.

Now we prove that in the junction tree over $X$ the number of clusters which
contain a variable $X_{i}$ is greater with 1 than the number of separators
which contain the same variable. This is true for all $i=1,\ldots,n$.
This means $\#\left\{  K\in\Gamma|X_{i}\in X_{K}\right\}  =$\#$\left\{
S\in\mathcal{S}|X_{i}\in X_{S}\right\}+1$.

For a variable $X_{i}$ we denote \#$\left\{ S\in\mathcal{S}|X_{i}\in
X_{S}\right\} $ by $t$.

Case: $t=0$.

The statement is a consequence of the definition of junction tree, that is
the union of all clusters is $X$, so every variable have to appear at least
in one cluster. $X_{i}$ can not appear in two clusters, because in this case
there should exist a separator which contain $X_{i}$ too, and we supposed
that there is not such a separator $(t=0)$

Case: $t>0$

If two clusters contain the variable $X_{i}$, then every cluster from the
path between the two clusters contain $X_{i}$ (running intersection
property). From this results that the clusters containing $X_{i}$ are the
nodes of a connected graph, and this graph is a tree. If this tree contain $%
t $ separator sets then it contains $t+1$ clusters. All of these separators
contain $X_{i}$, and each separator connects two clusters. So there will be $%
t+1$ clusters that contain $X_{i}$.

Applying this result in formula (\ref{eq:4}) after simplification we obtain

\[
f_{\mathbf{X}}\left(  \mathbf{x}\right)  =\dfrac{\prod\limits_{K\in\Gamma
}c_{\mathbf{X}_{K}}\left(  \mathbf{u}_{K}\right)  \prod\limits_{i=1}%
^{n}f_{X_{i}}\left(  x_{i}\right)  }{\prod\limits_{S\in\mathcal{S}}\left[
c_{\mathbf{X}_{S}}\left(  \mathbf{u}_{S}\right)  \right]  ^{v_{S}-1}}.
\]

Dividing both sides by $\prod\limits_{i=1}^{n}f_{X_{i}}\left(  x_{i}\right)$ we obtain

\begin{equation}\label{eq:5}
\dfrac{f_{\mathbf{X}}\left(  \mathbf{x}\right)  }{\prod\limits_{i=1}%
^{n}f_{X_{i}}\left(  x_{i}\right)  }=\dfrac{\prod\limits_{K\in\Gamma
}c_{\mathbf{X}_{K}}\left(  \mathbf{u}_{K}\right)  }{\prod\limits_{S\in
\mathcal{S}}\left[  c_{\mathbf{X}_{S}}\left(  \mathbf{u}_{S}\right)  \right]
^{v_{S}-1}}.
\end{equation}

Equations (\ref{eq:2}) and  (\ref{eq:5}) prove the statement of the theorem.

\end{pf}

We saw that if the conditional independence structure underlying the random
variables makes possible the construction of a junction tree, then the
multivariate copula density associated to the joint probability distribution
can be expressed as a product and fraction of lower dimensional copula
densities.

A logical question is the following. 
What conditions are necessary for (\ref{eq:5}) to be a copula density?
It is easy to see that
the product and fraction of copulas are positive. So 
\[
c\left( \mathbf{u}\right) =\dfrac{\prod\limits_{K\in \Gamma }c_{_K}\left( \mathbf{u}_K\right)
}{\prod\limits_{S\in \mathcal{S}}\left[ c_{_S}\left( \mathbf{u}_S\right) %
\right] ^{v_S-1}}
\] 
will be a copula density if and only if 
\[
\int\limits_{\left[ 0;1\right] ^n}\dfrac{\prod\limits_{K\in \Gamma
}c_{_K}\left( \mathbf{u}_K\right) }{\prod\limits_{S\in \mathcal{S}}\left[
c_{_S}\left( \mathbf{u}_S\right) \right] ^{v_S-1}} d\mathbf{u}=1.
\] 
This happens if the following consistency conditions are fulfilled for all connected
clique pairs $K_i$ and $K_j$: 

\[
\int\limits_{\left[ 0;1\right] ^{\#\left\{ K_{i}\backslash S_{ij}\right\}
}}c_{_{Ki}}\left( \mathbf{u}_{K_{i}}\right) d\mathbf{u}_{K_{i}\backslash
S_{ij}}=\int\limits_{\left[ 0;1\right] ^{\#\left\{ K_{j}\backslash
S_{ij}\right\} }}c_{_{Kj}}\left( \mathbf{u}_{K_{j}}\right) d\mathbf{u}%
_{K_{j}\backslash S_{ij}},
\]
where $S_{ij}=K_i$ $\cap $ $K_j$.
We emphasize
here that all cliques are subsets of the set $\Gamma$ of the acyclic
hypergraph $\left( V,\Gamma ,\mathcal{S}\right)$.

These conditions are fulfilled if $c_{_{S_{ij}}}\left( \mathbf{u}%
_{S_{ij}}\right) $ are marginal probability densities of $c_{K_{i}}\left(
\mathbf{u}_{K_{i}}\right)$, whenever $S_{ij}$
connects a cluster $K_{i}$. This can be expressed by terms of copula
function as follows.

For $\left\{ k_{1},\ldots ,k_{m}\right\} =K_{i}$ and $\left\{ s_{1},\ldots
,s_{l}\right\} =S_{ij}$, $\left\{ s_{1},\ldots ,s_{l}\right\} \subset
\left\{ k_{1},\ldots ,k_{m}\right\} $ stands $C_{m}\left( u_{k_{1}},\ldots
,u_{k_{m}}\right) =C_{l}\left( u_{s_{1}},\ldots ,u_{s_{l}}\right) $ for $%
u_{k_{i}}=1$, when $k_{i}\notin S_{ij}$. Usually this condition is not
fulfilled by copulas.

Finding multivariate copulas which
fulfill the consistency conditions is not a trivial task.

A special type of conditional independence, when the graph underlying the
random variables is starlike, is treated in {\cite{YaQiRu09}. Another 
type of special multivariate copula where the underlying conditional
independence graph is a tree can be found in \cite{Ki07}. 

For discrete random variables, the conditional independences are exploited
by the Markov random fields. In physics for two-valued random variables it
is known the Ising model. In these cases the random variables take on a few
values only. However many times the problem is hard. The great advantage
of using the discrete approach is that the marginal probability
distributions involved fulfill the consistency conditions ( see \cite{Vo62} and \cite{Vo63}).

If we have an i.i.d. sample of size $N$ from a continuous joint
probability distribution then for each random variable we have $N$
different values. For this
case, the empirical copulas were introduced and first studied by P. Deheuvels in \cite{De79} 
who called them empirical dependence functions. Later in \cite{Ma05} and \cite{Me05} there were introduced 
the so called discrete copulas. About the
two-dimensional empirical copulas one can read in Nelsen's introductory book (see \cite{Ne99}). 
In the case when we are dealing with a sample drawn from a continuous joint probability
distribution the size of these
random variables would be too large, so we will apply a uniform partition  
and define the so called sample derivated copula.

\section{The sample derivated copula.}
\label{sec:3}

Let $X_{1},\ldots ,X_{n}$ be continuous random variables in the
same probability field. 
Let
\begin{equation}\label{eq:6}
\begin{array}{c}
x_{1}^{1},\ldots ,x_{n}^{1}\\
x_{1}^{2},\ldots ,x_{n}^{2}\\
\vdots \\
x_{1}^{N},\ldots ,x_{n}^{N}\\
\end{array}
\end{equation}
be an i.i.d. sample of size $N$ taken from the joint probability distribution 
of the random vector $\left( X_{1},\ldots ,X_{n}\right)^{T}$.

As any sample element occurs two times in the sample with probability zero,
we can suppose that the sample elements are different.

We denote the set of the values of $X_{i}$ in the sample by $\Lambda _{i}$. This set
contains $N$ values, for each random variable. The theoretical range of the
continuous random variable $X_{i}$ will be denoted by $\overline{\Lambda _{i}}$. For
every $i$ we denote by $\lambda _{i}^{m}=\min \overline{\Lambda _{i}}\in R$
and by $\lambda _{i}^{M}=\max \overline{\Lambda _{i}}\in R$. We suppose for
simplicity that $\min \overline{\Lambda _{i}}\neq $ $\min \Lambda _{i}$ and $%
\max \overline{\Lambda _{i}}\neq \max \Lambda _{i}$ For each random variable
$X_{i}$ we define a partition of $\Lambda _{i}$ by $\mathcal{P}_{i}=\left\{ x_{0}^{p_{i}}=\lambda
_{i}^{m},x_{1}^{p_{i}},\ldots ,x_{m_{i}-1}^{p_{i}},x_{m_{i}}^{p_{i}}=\lambda
_{i}^{M}\right\} $ with the following properties:

\begin{itemize}
\item  For each random variable $X_{i}$, each interval $\left(
x_{j-1}^{p_{i}};x_{j}^{p_{i}}\right] ,j=1,\ldots, m_{i}$ contains a given $%
n_{i}=\dfrac{N}{m_{i}}\in N$ number of values from the set $\Lambda _{i}$.

\item  Each $x_{j}^{p_{i}}\in \Lambda _{i},\ j=1,\ldots, m_{i}-1$.
\end{itemize}

The  partition with the above properties will be called uniform partition.
We denote by ${\mathcal P}$ the set of partitions $\{{\mathcal P}_{1},\ldots, {\mathcal P}_{n}\}$.

Let be $\widetilde{X_{i}}$ the categorical random variable associated to the
random variable $X_{i}$:

\[
P(\widetilde{X_{i}} \in \left( x_{j-1}^{p_{i}};x_{j}^{p_{i}}\right])=\dfrac{1}{m_{i}}, j=1,\ldots , m_{i}.
\]

We assign to each $x^{i}\in \left( x_{j-1}^{p_{i}};x_{j}^{p_{i}}\right] $
the number $u_{j}^{i}=\dfrac{j}{m_{i}}, j=0,\ldots, m_{i}$. Obviously
$u_{0}^{i}=0$ and $u_{m_{i}}^{i}=1$. Let 
$\widetilde{\Lambda }_{i}=\left\{ u_{j}^{i}|j=0,\ldots , m_{i}\right\} $.
So we can define the following discrete uniform random variables:

\[
\widetilde{U_{i}}=
\left(
\begin{array}{ccccc}
u_{0}^{i} & u_{1}^{i}      & \ldots & u_{m_{i}-1}^{i} & u_{m_{i}}^{i}   \vspace{3mm}\\
0         &\dfrac{1}{m_{i}} & \ldots & \dfrac{1}{m_{i}} & \dfrac{1}{m_{i}}
\end{array}
\right),
i=1,\ldots, n.
\] 
Now we transform the sample (\ref{eq:6}) using the above assignment. 
We denote the transformed sample by ${\mathcal T}$. 

\begin{dfn}\label{de:4} 
The function 
$\widetilde{c}:\prod\limits_{i=1}^{n}\widetilde{\Lambda }_{i}\rightarrow R$ 
defined by
\[
\left( u_{k_{1}}^{1},\ldots ,u_{k_{n}}^{n}\right) \mapsto 
\widetilde{c}\left(
u_{k_{1}}^{1},\ldots ,u_{k_{n}}^{n}\right) =\dfrac{\#\left\{(u_{k_{1}}^{1},\ldots ,u_{k_{n}}^{n})\in {\mathcal T}\right\}}{N},
k_{i}=0,\ldots, m_{i}  
\]
will be called sample derivated copula distribution.
\end{dfn}

\begin{rmk}\label{re:2}
The maximum number of different vectors what the above defined sample derivated copula can take on equals to 
$\prod\limits_{i=1}^{n}m_{i}$.
\end{rmk}

\begin{dfn}\label{de:5} 
The function 
$\widetilde{C}^{\mathcal{P}}_{n}:
\prod\limits_{i=1}^{n}\widetilde{\Lambda }_{i}\subset \left[
0;1\right] ^{n}\rightarrow \left[ 0;1\right]$
defined by
\[
\begin{array}{l}
\left( u_{k_{1}}^{1},\ldots ,u_{k_{n}}^{n}\right) \mapsto 
\widetilde{C}^{\mathcal{P}}_{n}\left(
u_{k_{1}}^{1},\ldots, u_{k_{n}}^{n}\right) = \vspace{3mm}\\
= \dfrac{\#\left\{(u_{k_{1}}^{1},\ldots ,u_{k_{n}}^{n})\in {\mathcal T} \;|\; u_{1} \le u_{k_{1}}^{1},\ldots ,u_{n} \le u_{k_{n}}^{n}\right\}}{N}
\end{array}
\]
will be called sample derivated copula.
\end{dfn}

Throughout the paper we use the notation $\widetilde{C}_{n}$ instead of $\widetilde{%
C}^{\mathcal{P}}_{n}$.

\begin{thm}\label{theo:2} 
The sample derivated copula is a copula.
\end{thm}

\begin{pf}

\begin{itemize}

\item[1)] It is evident that $\widetilde{C}_{n}$ is increasing in its each component.

\item[2)] If exists $s$ such that $u_{k_{s}}^{s}=0$ then $\widetilde{C}_{n}\left(
u_{k_{1}}^{1},\ldots ,u_{k_{s-1}}^{s-1},0,u_{k_{s+1}}^{s+1}\ldots
,u_{k_{n}}^{n}\right) =0$. This follows directly from the definition. The
sample do not contain any vector with a negative coordinate.

\item[3)] If for all $s\neq l$ we have $u_{k_{s}}^{s}$ $=1$ then

\[
\begin{array}{l}
\widetilde{C}_{n}\left( 1,\ldots 1,u_{k_{l}}^{l},1\ldots ,1\right) = \\
 = \dfrac{\#\left\{ \left( u_{1},\ldots ,u_{s},\ldots u_{n}\right) \in {\mathcal T} \; | \; u_{s}\le1, \forall s\neq l\text{ and }%
u_{l} \le u_{k_{l}}^{l}\right\} }{N}=\\
 = \dfrac{1}{N}\sum\limits_{0}^{k_{l}-1}n_{l}=\dfrac{n_{l}\cdot k_{l}}{N}=\dfrac{%
k_{l}}{m_{l}}=u_{k_{l}}^{l}
\end{array}
\]

\item[4)] $\widetilde{C}_{n}$ is $n$-increasing as it is a
cumulative probability distribution function.

\end{itemize}

\end{pf}

\begin{rmk}\label{re:3}
The sample derivated copula differs, from the empirical copula
\cite{De79} and the discrete copula \cite{Ma05}, \cite{Me05}. One of the differences is that the
cardinal of $\widetilde{\Lambda }_{i}$ is not necessary the same for all $i=1,\ldots, n$. 
Another difference is that a marginal variable can take the same value in
more than one vector (since $m_{i}<N$).
\end{rmk}

\begin{thm}\label{theo:3} 
The sample derivated copula has the following consistency property.
If all variables $u_{k_{s}}^{s}$ $=1$ for $s\in V\backslash \left\{ l_{1},\ldots ,l_{q}\right\}$ 
then

\[
\widetilde{C}_{n}\left( u_{k_{1}}^{1},\ldots, u_{k_{n}}^{n}\right) 
=\widetilde{C}_{q}\left( u_{l_{1}}^{1},\ldots, u_{l_{q}}^{q}\right).
\]
\end{thm}

\begin{pf}
\[
\begin{array}{l}
\widetilde{C}_{n}\left( u_{k_{1}}^{1},\ldots,u_{k_{n}}^{n}\right) = \vspace{3mm}\\
= \dfrac{\#\left\{ \left( u_{1},\ldots, u_{n}\right) \in {\mathcal T}\;|\;u_{s}\le 1, s\in V
\backslash \left\{ l_{1},\ldots ,l_{q}\right\},
u_{l_{i}}\le u_{l_{i}}^{i},i=1\ldots q\right\} }{N}=\vspace{3mm}\\
=\dfrac{\#\left\{ \left( u_{1},\ldots, u_{n}\right) \in {\mathcal T} \;|\;
u_{l_{i}}\le u_{l_{i}}^{i},i=1\ldots q\right\} }{N}=\vspace{3mm}\\
=\widetilde{C}_{q}\left( u_{l_{1}}^{1},\ldots ,u_{l_{q}}^{q}\right).
\end{array}
\]

\end{pf}

\begin{rmk}\label{re:4}
In general copulas do not fulfill the consistency property..
\end{rmk}
\begin{rmk}\label{re:5}
This theorem assures the consistency statements that we need when
constructing junction tree like copulas.
\end{rmk}

At the end of this part we convince the reader from an information
theoretical point of view why should one use the uniform partition and
the sample derivated copula.

In the following we suppose that each $\Lambda _{i}, i=1,\ldots, n$ is partitioned in the
same number of $m_{i}, i=1,\ldots, n$ intervals as in the previous case. We denote now the partitioning 
points by $y_{j}^{p_{i}},\; j=0, 1, \ldots, m_{i}; \; i=1, \ldots, n$. 
This partition is arbitrary (for example equidistant) which has not the property that each interval
contains the same number of sample elements. The partition of $\Lambda
_{i}$ is given by $\left\{ y_{j}^{p_{i}}\;|\;j\in \left\{ 0, 1, \ldots, m_{i}\right\} \right\}$ 
and is denoted by $\mathcal{P}_{i}^{\prime }$. We denote by ${\mathcal P}^{\prime}$ the set of
partitions ${\mathcal P}_{1}^{\prime}, \ldots, {\mathcal P}_{n}^{\prime}$.

We denote the number of values of the variable $X_{i}$ belonging to $\left(
y_{j}^{p_{i}};y_{j+1}^{p_{i}}\right] \cap \Lambda _{i}$ by $k_{j}^{i}$.

Let $\widetilde{Y}_{i}$ be the categorical random variable associated to 
$X_{i}$:
\[
P\left( \widetilde{Y}_{i}\in \left( y_{j}^{p_{i}};y_{j+1}^{p_{i}}\right]
\right) =\dfrac{k_{j}^{i}}{N}, j=0, 1, \ldots, m_{i},
\]
where
\[
\sum\limits_{j=1}^{m_{i}}\dfrac{k_{j}^{i}}{N}=1.
\]
The entropy of $X_{i}$ determined by the partition 
$\mathcal{P}_{i}^{\prime }$ is: 
\[
H_{\mathcal{P}_{i}^{\prime }}\left( X_{i}\right)
=H\left( \widetilde{Y}_{i}\right) = - \sum\limits_{j=1}^{m_{i}}\dfrac{k_{j}^{i}}{N}%
\log \dfrac{k_{j}^{i}}{N},\; i=1\ldots n.
\]
It can be seen that the entropy $H_{\mathcal{P}_{i}^{\prime }}\left(X_{i}\right)$
depends on the number of intervals $m_{i}$ and on $k_{j}^{i}$.

We introduce the following notation:
\[
q_{j_{1},\ldots, j_{n}}^{X_{1},\ldots
,X_{n}}=P\left( X_{1}\in \left( y_{j_{1}-1}^{p_{1}};y_{j_{1}}^{p_{1}}\right]
,X_{2}\in \left( y_{j_{2}-1}^{p_{2}};y_{j_{2}}^{p_{2}}\right] ,\ldots
,X_{n}\in \left( y_{j_{n}-1}^{p_{n}};y_{j_{n}}^{p_{n}}\right] \right),
\]
where $j_{i}=1,\ldots,m_{i},\; i=1,\ldots,n.$

The joint probability distribution determined by the partition 
$\mathcal{P}^{^{\prime }}$ has the joint entropy:
\[
H_{ \mathcal{P}^{\prime }}\left(
X_{1},\ldots ,X_{n}\right) =-\sum\limits_{j_{1}=1}^{m_{1}}\cdots\sum\limits_{j_{n}=1}^{m_{n}}q_{j_{1},\ldots j_{n}}^{X_{1},\ldots
,X_{n}}\log _{2}q_{j_{1},\ldots j_{n}}^{X_{1},\ldots ,X_{n}}.
\]
The information content of the joint probability distribution determined by the partition 
$\mathcal{P}^{^{\prime }}$  is:
\begin{equation}\label{eq:7}
\begin{array}{l}
I_{ \mathcal{P}^{\prime }}\left(
X_{1},\ldots ,X_{n}\right) =\sum\limits_{i=1}^{n}H_{\mathcal{P}_{i}^{\prime
}}\left( X_{i}\right) -H_{ \mathcal{P}^{\prime }}\left( X_{1},\ldots ,X_{n}\right)=\vspace{3mm}\\
=-\sum\limits_{i=1}^{n}\sum\limits_{j=1}^{m_{i}}\dfrac{k_{j}^{i}}{N}\log 
\dfrac{k_{j}^{i}}{m_{i}}
+\sum\limits_{j_{1}=1}^{m_{1}}\cdots\sum\limits_{j_{n}=1}^{m_{n}}q_{j_{1},\ldots j_{n}}^{X_{1},\ldots ,X_{n}}\log
_{2}q_{j_{1},\ldots j_{n}}^{X_{1},\ldots ,X_{n}}.
\end{array}
\end{equation}

\begin{rmk}\label{re:6}
In this case the information content depends on the number of
intervals $m_{i}$, the number of values in each one dimensional interval $k_{j}^{i}$, and the probabilities of 
belonging to the $n$-dimensional intervals.
\end{rmk}

If we regard the uniform partition $\mathcal{P}$ then the entropy of $X_{i}$ is:
\[
H_{\mathcal{P}_{i}}\left( X_{i}\right) =H\left( \widetilde{X}_{i}\right)
=-\sum\limits_{j=1}^{m_{i}}\dfrac{1}{m_{i}}\log \dfrac{1}{m_{i}}=\log m_{i}, \; i=1,\ldots ,n.
\]
The entropy $H_{\mathcal{P}_{i}}\left( X_{i}\right) $ depends just on the
number of intervals $m_{i}$.

We express now the probability
\[
\begin{array}{l}
p_{j_{1},\ldots j_{n}}^{X_{1},\ldots ,X_{n}}=P\left( X_{1}\in \left(
x_{j_{1}-1}^{p_{1}};x_{j_{1}}^{p_{1}}\right] ,\ldots ,X_{n}\in \left(
x_{j_{n}-1}^{p_{n}};x_{j_{n}}^{p_{n}}\right] \right) =\vspace{3mm}\\
P\left( \widetilde{U}_{1}=u_{j_{1}},\ldots
,\widetilde{U}_{n}=u_{j_{n}}\right) =\widetilde{c}\left( u%
_{j_{1}},\ldots ,u_{jn}\right) 
\end{array}
\]

The joint probability entropy associated to the partition is:
\[
\begin{array}{l}
H_{\mathcal{P}}\left( X_{1},\ldots
,X_{n}\right) =-\sum\limits_{j_{1}=1}^{m_{1}}\cdots\sum\limits_{j_{n}=1}^{m_{n}}
p_{j_{1},\ldots j_{n}}^{X_{1},\ldots ,X_{n}}\log
_{2}p_{j_{1},\ldots j_{n}}^{X_{1},\ldots ,X_{n}}=\vspace{3mm}\\
=-\sum\limits_{j_{1}=1}^{m_{1}}\cdots\sum\limits_{j_{n}=1}^{m_{n}}
\widetilde{c}\left( u_{j_{1}},\ldots ,u_{jn}\right) \log _{2}\widetilde{c}\left( 
u_{j_{1}},\ldots ,u%
_{jn}\right).
\end{array}
\]

The information content determined by the partition $\mathcal{P}$ is:
\begin{equation}\label{eq:8}
\begin{array}{l}
I_{\mathcal{P}}\left(
X_{1},\ldots ,X_{n}\right) =\sum\limits_{i=1}^{n}H_{\mathcal{P}_{i}}\left( X_{i}\right) 
- H_{\mathcal{P}}\left( X_{1},\ldots ,X_{n}\right) =\vspace{3mm}\\
=\sum\limits_{i=1}^{n}\log m_{i}+\sum\limits_{j_{1}=1}^{m_{1}}\cdots\sum\limits_{j_{n}=1}^{m_{n}}
\widetilde{c}\left( u%
_{j_{1}},\ldots ,u_{jn}\right) \log _{2}%
\widetilde{c}\left( u_{j_{1}},\ldots ,%
u_{jn}\right)
\end{array}
\end{equation}

\begin{rmk}\label{re:7}
If we suppose that for all $i=1,\ldots, n$ the number of intervals 
$m_{i}$ is the same for the two discussed cases then comparing formulas (\ref{eq:7})
and (\ref{eq:8}) we can see that in the case of partition $\mathcal{P}$ the information content 
does not depend on the first sum of formula (\ref{eq:8}) but only on the
sample derivated copula.
\end{rmk}

\section{The junction tree approach applied to the sample derivated copula.}
\label{sec:4}

We introduced the sample derivated copula as a discrete probability
distribution with uniform marginals. We proved for this special copula in Theorem \ref{theo:3} 
that the consistency properties are fulfilled.

Let $V=\left\{ 1,\ldots ,n\right\}$ be again a set of vertices. 
Let be defined an acyclic hypergraph over 
$V$. We denote by $\Gamma $ and $\mathcal{S}$ the set of clusters and
separators of the hypergraph which determine a junction tree ${\mathcal J}$. The marginal
probability distributions associated to the clusters 
$K=\left\{ i_{1},\ldots ,i_{t}\right\}\in \Gamma $ are
denoted by $\widetilde{c}_{K}\left( \widetilde{\mathbf{U}}_{K}\right) =$ $%
\widetilde{c}_{i_{1},\ldots ,i_{t}}\left( \widetilde{U}_{i_{1}},\ldots ,%
\widetilde{U}_{i_{t}}\right)$. The marginal probability distributions
associated to the separators are denoted in the same way by 
$\widetilde{c}_{S}\left( \widetilde{\mathbf{U}}_{S}\right) $. 
The joint discrete copula is shortly denoted by $\widetilde{%
c}\left( \widetilde{\mathbf{U}}\right)$ and the univariate marginals by $%
\widetilde{c}_{i}\left( \widetilde{U}_{i}\right) ,i=1,\ldots ,n$ . 

In this section we are going to use the following popular notation:
\[
\sum\limits_{\mathbf{u}}f\left( \widetilde{\mathbf{U}}\right)
=\sum\limits_{i_{1}=1}^{m_{1}}\ldots \sum\limits_{i_{n}=1}^{m_{n}}f\left( 
\widetilde{U}_{1}=u_{i_{1}}^{1},\ldots ,\widetilde{U}_{n}=u_{i_{n}}^{n}%
\right),
\]
where $u_{i_{k}}^{k},i_{k}=1,\ldots ,m_{k}$ are the possible values of the
random variable $\widetilde{U}_{k}, k=1, \ldots, n$ and $f$ is an arbitrary 
$n$-dimensional function. This simplified notation is used for the products, too.

\begin{dfn}\label{de:6}
The junction tree distribution given by 
\begin{equation}\label{eq:9}
\widetilde{c}_{\mathcal J}\left( \widetilde{\mathbf{U}}\right) =\dfrac{%
\prod\limits_{K\in \Gamma }\widetilde{c}_{K}\left( \widetilde{\mathbf{U}}%
_{K}\right) }{\prod\limits_{S\in \mathcal{S}}\left[ \widetilde{c}_{S}\left( 
\widetilde{\mathbf{U}}_{S}\right) \right] ^{v_{S}-1}},
\end{equation}
where $v_{S}$ is
the number of clusters connected by the separator $S$,
is called copula junction tree distribution, or shortly junction tree copula.
\end{dfn}

The problem is finding the junction tree copula which fits to the
sample derivated copula. The goodness of fitting will be quantified by the
Kullback-Leibler divergence \cite{CoTo91}.

\begin{thm}\label{theo:4} 
The Kullback-Leibler divergence between the approximation (\ref{eq:9}) and
the sample derivated copula $\widetilde{c}\left( \widetilde{\mathbf{U}}\right) $%
is given by the formula:
\[
\begin{array}{l}
KL\left( \widetilde{c}_{\mathcal J}\left( \widetilde{\mathbf{U}}%
\right) ,\widetilde{c}\left( \widetilde{\mathbf{U}}\right) \right) =\vspace{3mm}\\
= -H\left( 
\widetilde{\mathbf{U}}\right) -\left[ \sum\limits_{K\in \Gamma }I\left( 
\widetilde{\mathbf{U}}_{K}\right) -\sum\limits_{S\in \mathcal{S}}\left(
v_{S}-1\right) I\left( \widetilde{\mathbf{U}}_{S}\right) \right]
+\sum\limits_{i=1}^{n}\log _{2}m_{i}.
\end{array}
\]
\end{thm}

\begin{pf}

\[
\begin{array}{l}
KL\left( \widetilde{c}_{\mathcal J}\left( \widetilde{\mathbf{U}}\right) ,%
\widetilde{c}\left( \widetilde{\mathbf{U}}\right) \right) =\vspace{3mm}\\
=\sum\limits_{%
\mathbf{u}}\widetilde{c}\left( \widetilde{\mathbf{U}}\right) \log _{2}\dfrac{%
\widetilde{c}\left( \widetilde{\mathbf{U}}\right) }{\widetilde{c}_{\Gamma
}\left( \widetilde{\mathbf{U}}\right) }=\sum\limits_{\mathbf{u}}\widetilde{c}%
\left( \widetilde{\mathbf{U}}\right) \log _{2}\widetilde{c}\left( \widetilde{%
\mathbf{U}}\right) -\sum\limits_{\mathbf{u}}\widetilde{c}\left( \widetilde{%
\mathbf{U}}\right) \log _{2}\widetilde{c}_{\Gamma }\left( \widetilde{\mathbf{%
U}}\right) =\vspace{3mm}\\
=-H\left( \widetilde{\mathbf{U}}\right) -\sum\limits_{\mathbf{u}}\widetilde{%
c}\left( \widetilde{\mathbf{U}}\right) \log _{2}\dfrac{\prod\limits_{K\in
\Gamma }\widetilde{c}_{K}\left( \widetilde{\mathbf{U}}_{K}\right) }{%
\prod\limits_{S\in \mathcal{S}}\left[ \widetilde{c}_{S}\left( \widetilde{%
\mathbf{U}}_{S}\right) \right] ^{v_{S}-1}}=-H\left( \widetilde{\mathbf{U}}%
\right) - \vspace{3mm}\\
-\sum\limits_{\mathbf{u}}\widetilde{c}\left( \widetilde{\mathbf{U}}\right) %
\left[ \log _{2}\prod\limits_{K\in \Gamma }\widetilde{c}_{K}\left( 
\widetilde{\mathbf{U}}_{K}\right) -\log _{2}\prod\limits_{S\in \mathcal{S}}%
\left[ \widetilde{c}_{S}\left( \widetilde{\mathbf{U}}_{S}\right) \right]
^{v_{S}-1}\right] = \vspace{3mm}\\
=-H\left( \widetilde{\mathbf{U}}\right) -\sum\limits_{\mathbf{u}}\widetilde{%
c}\left( \widetilde{\mathbf{U}}\right) \log _{2}\prod\limits_{K\in \Gamma }%
\widetilde{c}_{K}\left( \widetilde{\mathbf{U}}_{K}\right) +\vspace{3mm}\\
+\sum\limits_{%
\mathbf{u}}\widetilde{c}\left( \widetilde{\mathbf{U}}\right) \log
_{2}\prod\limits_{S\in \mathcal{S}}\left[ \widetilde{c}_{S}\left( \widetilde{%
\mathbf{U}}_{S}\right) \right] ^{v_{S}-1}.
\end{array}
\]
We add and substract the sum:
\begin{equation}\label{eq:10}
\sum\limits_{\mathbf{u}}\widetilde{c}\left( \widetilde{\mathbf{U}}\right)
\log _{2}\prod\limits_{K\in \Gamma }\prod\limits_{i\in K}\widetilde{c}%
_{i}\left( \widetilde{U}_{i}\right).
\end{equation}
It follows from the definition of the junction tree that 
$\bigcup\limits_{K\in \Gamma }^{{}}K=V$, and each variable belongs once more
in the clusters as in the separators. So by adding and substracting (\ref{eq:10}) we
obtain the following:
\[
\begin{array}{l}
KL\left( \widetilde{c}_{\mathcal J}\left( \widetilde{\mathbf{U}}\right) ,%
\widetilde{c}\left( \widetilde{\mathbf{U}}\right) \right) =-H\left( 
\widetilde{\mathbf{U}}\right) -\sum\limits_{\mathbf{u}}\widetilde{c}\left( 
\widetilde{\mathbf{U}}\right) \log _{2}\dfrac{\prod\limits_{K\in \Gamma }%
\widetilde{c}_{K}\left( \widetilde{\mathbf{U}}_{K}\right) }{%
\prod\limits_{K\in \Gamma }\prod\limits_{i\in K}\widetilde{c}_{i}\left( 
\widetilde{U}_{i}\right) }+\vspace{3mm}\\
+\sum\limits_{\mathbf{u}}\widetilde{c}\left( \widetilde{\mathbf{U}}\right)
\log _{2}\dfrac{\prod\limits_{S\in \mathcal{S}}\left[ \widetilde{c}_{S}\left( 
\widetilde{\mathbf{U}}_{S}\right) \right] ^{v_{S}-1}}{\prod\limits_{S\in 
\mathcal{S}}\left[ \prod\limits_{i\in S}\widetilde{c}_{i}\left( \widetilde{U}%
_{i}\right) \right] ^{v_{S}-1}}-\sum\limits_{\mathbf{u}}\widetilde{c}\left( 
\widetilde{\mathbf{U}}\right) \log _{2}\prod\limits_{i=1}^{n}\widetilde{c}%
_{i}\left( \widetilde{U}_{i}\right) =
\end{array}
\]
\[
\begin{array}{l}
=-H\left( \widetilde{\mathbf{U}}\right) -\sum\limits_{\mathbf{u}}\widetilde{%
c}\left( \widetilde{\mathbf{U}}\right) \sum\limits_{K\in \Gamma }\log _{2}%
\dfrac{\widetilde{c}_{K}\left( \widetilde{\mathbf{U}}_{K}\right) }{%
\prod\limits_{i\in K}\widetilde{c}_{i}\left( \widetilde{U}_{i}\right) }+\vspace{3mm}\\
+\sum\limits_{\mathbf{u}} \widetilde{c}\left( \widetilde{\mathbf{U}}\right) \sum\limits_{S\in 
\mathcal{S}}\log _{2}\dfrac{\left[ \widetilde{c}_{S}\left( \widetilde{\mathbf{%
U}}_{S}\right) \right] ^{v_{S}-1}}{\left[ \prod\limits_{i\in S}\widetilde{c}%
_{i}\left( \widetilde{U}_{i}\right) \right] ^{v_{S}-1}}-\sum\limits_{\mathbf{%
u}}\widetilde{c}\left( \widetilde{\mathbf{U}}\right)
\sum\limits_{i=1}^{n}\log _{2}\widetilde{c}_{i}\left( \widetilde{U}%
_{i}\right).
\end{array}
\]

Since the sample derivated copula has the property that all $\widetilde{c}%
_{K}\left( \widetilde{\mathbf{U}}_{K}\right)$, $\widetilde{c}_{S}\left( 
\widetilde{\mathbf{U}}_{S}\right)$, $\widetilde{c}_{i}\left( \widetilde{U}%
_{i}\right)$ are consistent marginals of $\widetilde{c}\left( \widetilde{\mathbf{U}}%
\right) $ (see Theorem \ref{theo:3}) we have the following relations:

\[
\begin{array}{l}
\bullet\;\;
\sum\limits_{\mathbf{u}}\widetilde{c}\left( \widetilde{\mathbf{U}}%
\right) \sum\limits_{K\in \Gamma }\log _{2}\dfrac{\widetilde{c}_{K}\left( 
\widetilde{\mathbf{U}}_{K}\right) }{\prod\limits_{i\in K}\widetilde{c}%
_{i}\left( \widetilde{U}_{i}\right) }=\vspace{3mm}\\
=\sum\limits_{K\in \Gamma }\sum\limits_{%
\mathbf{u}_{k}}\widetilde{c}_{K}\left( \widetilde{\mathbf{U}}_{K}\right)
\log _{2}\dfrac{\widetilde{c}_{K}\left( \widetilde{\mathbf{U}}_{K}\right) }{%
\prod\limits_{i\in K}\widetilde{c}_{i}\left( \widetilde{U}_{i}\right) }%
=\sum\limits_{K\in \Gamma }I\left( \widetilde{\mathbf{U}}_{K}\right);\vspace{3mm}
\\
\bullet\;\;
\sum\limits_{\mathbf{u}}\widetilde{c}\left( \widetilde{\mathbf{U}}%
\right) \sum\limits_{S\in \mathcal{S}}\log _{2}\dfrac{\left[ \widetilde{c}%
_{S}\left( \widetilde{\mathbf{U}}_{S}\right) \right] ^{v_{S}-1}}{\left[
\prod\limits_{i\in S}\widetilde{c}_{i}\left( \widetilde{U}_{i}\right) \right]
^{v_{S}-1}}=\vspace{3mm}\\
=\sum\limits_{S\in \mathcal{S}}\sum\limits_{\mathbf{u}_{S}}\left(
v_{S}-1\right) \widetilde{c}_{S}\left( \widetilde{\mathbf{U}}_{S}\right)
\log _{2}\dfrac{\widetilde{c}_{S}\left( \widetilde{\mathbf{U}}_{S}\right) }{%
\prod\limits_{i\in S}\widetilde{c}_{i}\left( \widetilde{U}_{i}\right) }
=\sum\limits_{S\in \mathcal{S}}\left( v_{S}-1\right) I\left( \widetilde{%
\mathbf{U}}_{S}\right);\vspace{3mm}
\\
\bullet\;\;
-\sum\limits_{\mathbf{u}}\widetilde{c}\left( \widetilde{\mathbf{U}}%
\right) \sum\limits_{i=1}^{n}\log _{2}\widetilde{c}_{i}\left( \widetilde{U}%
_{i}\right) =\sum\limits_{i=1}^{n}H\left( \widetilde{U}_{i}\right)=\sum\limits_{i=1}^{n}\log _{2}m_{i};
\end{array}
\]
Here $I\left( \widetilde{\mathbf{U}}_{K}\right) ,I\left( \widetilde{\mathbf{%
U}}_{S}\right) $ are the information content of the probability distribution
of the marginals $\widetilde{c}_{K}\left( \widetilde{\mathbf{U}}_{K}\right) $
and $\widetilde{c}_{S}\left( \widetilde{\mathbf{U}}_{S}\right) $ (see \cite{CoTo91}).

By the substitution of these assertions we obtain:

\[
\begin{array}{l}
KL\left( \widetilde{c}_{\mathcal J}\left( \widetilde{\mathbf{U}}\right) ,%
\widetilde{c}\left( \widetilde{\mathbf{U}}\right) \right) = \vspace{3mm}\\
= -H\left( 
\widetilde{\mathbf{U}}\right) -\left[ \sum\limits_{K\in \Gamma }I\left( 
\widetilde{\mathbf{U}}_{K}\right) -\sum\limits_{S\in \mathcal{S}}\left(
v_{S}-1\right) I\left( \widetilde{\mathbf{U}}_{S}\right) \right]
+\sum\limits_{i=1}^{n}\log _{2}m_{i}.
\end{array}
\]

\end{pf}

\begin{rmk}\label{re:8}
The difference $\sum\limits_{i=1}^{n}\log _{2}m_{i}-H\left( 
\widetilde{\mathbf{U}}\right) $ does not depend on the junction tree
structure.
\end{rmk}

\begin{dfn}\label{de:7}
The difference 
\[
\sum\limits_{K\in \Gamma }I\left( \widetilde{\mathbf{U}}_{K}\right)
-\sum\limits_{S\in \mathcal{S}}\left( v_{S}-1\right) I\left( \widetilde{%
\mathbf{U}}_{S}\right)
\]
is called the weight of the junction tree copula.
\end{dfn}

It is easy to see that in order to find a better approximation using
junction trees, wee have to find the junction tree having the largest
weight. 

Finding the best fitting $k$-width junction tree, (the largest cluster
contains $k$ elements) for $k > 2$ is an NP-hard problem. For $k=2$
the problem is similar to the Chow-Liu approximation \cite{ChLi68}. In this case it is
possible to find the best fitting second order junction tree by Kruskal' or
Prim' algorithm.

For $k\geq 3$ it can be successfully used a heuristic approach introduced by
the authors in \cite{KoSza10} and \cite{SzaKo08}. The idea is the fitting of a special kind of
junction tree, called $t$-cherry junction tree. 

\section{Conclusions and possible applications}
\label{sec:5}

One of the advantages of the junction tree copula is that it reveals some of
the conditional independences between the variables involved. This kind of
dependence structure is not exploited by the copula function. Another
advantage of the method is that a multivariate copula can be decomposed into some
lower dimensional sample derivated copulas.

The sample derivated copula approach is useful in cases when nothing else is known about
the probability distribution but an iid sample. If the uniform partition is applied the whole
information content depends on the sample derivated copula. 

The copula junction tree can be used in feature selection which is a key-question
in many fields as finance, medicine and biostatistics.

We got very good numerical results in pattern recognition (see \cite{SzaKo10}). First we applied the uniform partition 
to discretize continuous random variables then constructed the $t$-cherry junction tree approximation.
In this way we found the informative features and so reduced the dimension of the classifier.

%% The Appendices part is started with the command \appendix;
%% appendix sections are then done as normal sections
%% \appendix

%% \section{}
%% \label{}

%% References
%%
%% Following citation commands can be used in the body text:
%% Usage of \cite is as follows:
%%   \cite{key}         ==>>  [#]
%%   \cite[chap. 2]{key} ==>> [#, chap. 2]
%%

%% References with bibTeX database:

%\bibliographystyle{elsarticle-num}
%\bibliography{<your-bib-database>}

%% Authors are advised to submit their bibtex database files. They are
%% requested to list a bibtex style file in the manuscript if they do
%% not want to use elsarticle-num.bst.

%% References without bibTeX database:

\end{document}